\documentclass[12pt, twoside, leqno]{article}

\usepackage{amsmath,amsthm}
\usepackage{amssymb}

\usepackage{enumerate}

\usepackage{graphicx}

\newtheorem{thm}{Theorem}[section]

\newtheorem{lem}[thm]{Lemma}

\theoremstyle{definition}
\newtheorem{defin}[thm]{Definition}

\numberwithin{equation}{section}

\usepackage{pgf,tikz}
\usetikzlibrary{arrows}

\begin{document}

\baselineskip=12.1pt


\title{ $L(n)$ graphs are vertex-pancyclic and Hamilton-connected}

\author{S.Morteza Mirafzal$^*$, Sara Kouhi\\
Department of Mathematics \\
  Lorestan University, Khorramabad, Iran\\
E-mail: mirafzal.m@lu.ac.ir\\
E-mail: smortezamirafzal@yahoo.com\\
\\Sara Kouhi\\Department of Mathematics \\
  Lorestan University, Khorramabad, Iran\\
E-mail:sarakouhi22@gmail.com
}

\date{}

\maketitle

\renewcommand{\thefootnote}{}


\footnote{\emph{*Corresponding author.}}
\footnote{\emph{Date}:  26 3 2022}
\renewcommand{\thefootnote}{\arabic{footnote}}

\setcounter{footnote}{0}


\begin{abstract} A graph $G$ of order $n>2$ is pancyclic if $G$ contains a cycle
of length $l$ for each integer $l$ with $3 \leq l \leq n $ and it is called vertex-pancyclic if every vertex is contained in a cycle of
length $l$ for every $3 \leq l \leq n $. A graph $G$ of order $n > 2$ is Hamilton-connected if for
any pair of distinct vertices $u$ and $v$, there is a Hamilton $u$-$v$ path, namely, there is
a $u$-$v$ path of length $n-1$.   The graph $ B(n)$ is a
 graph with the vertex set  $V=\{v \  |  \ v \subset [n] ,  | v |  \in \{ 1,2  \} \} $ and the
edge set $ E= \{ \{ v , w \} \  | \  v , w \in V , v \subset w $ or $ w \subset v \}$, where $[n]=\{1,2,...,n\}$. We denote by $L(n)$ the line graph of $B(n)$, that is, $L(n)=L(B(n))$.
  In this paper,  we show that the graph $L(n)$ is  vertex-pancyclic and Hamilton-connected  whenever $n\geq 6$.

\end{abstract}

\maketitle
\section{ Introduction and Preliminaries}  In this paper, a graph $G=(V,E)$ is
considered as an undirected simple graph where $V=V(G)$ is the vertex-set
and $E=E(G)$ is the edge-set. For all the terminology and notation
not defined here, we follow [3,4].\

 Let $n\geq 3$ be an integer and $[n]=\{1,2,...,n\}$.  The graph $ B(n)$ is a
 graph with the vertex set  $V=\{v \  |  \ v \subset [n] ,  | v |  \in \{ 1,2  \} \} $ and the
edge set $ E= \{ \{ v , w \} \  | \  v , w \in V , v \subset w $ or $ w \subset v \} $. We denote by $L(n)$ the line graph of $B(n)$, that is, $L(n)=L(B(n))$. It is easy to check that $L(n)$ is a connected regular graph of regularity $n-1$. Hence  when $n=3$,  $L(n)$ is the cycle $C_6$. It can easily be seen that the graph $B(n)$ is an edge-transitive graph [8], thus $L(n)$ is a
 vertex-transitive graph [3,4]. The graphs $B(n)$  and $L(n)$ have some interesting properties and have been studied in some aspects [6,7,8,9].  In [8] it has been shown that the graph $L(n)$ is a Cayley graph if and only if $n$  is a power of a prime. The graph $G$ is called an $integral$ graph whenever each of the eigenvalues of its adjacency matrix is an integer. In [6] it has been proved that the graph $L(n)$ is an integral graph. In fact the set of eigenvalues of $L(n)$ is $\{ -2,-1,0,n-2,n-1 \}$. An interesting property of the graph $B(n)$ has been appeared in [7], that is, $L(n)$ is  (isomorphic to)  the square root of the Johnson graph $J(n+1,2)$, namely,  $B(n)^2$ is isomorphic to $J(n+1,2)$.  If $v$ is a vertex of the graph $B(n)$, then $deg(v) \in \{2, n-1 \}$. Therefore, if $n$ is an odd integer then $B(n)$ is an $eulerian$  graph [4].  Hence $L(n)$ is a $Hamilton$  graph when $n$ is an odd integer. A graph $G$ of order $n > 2$ is $Hamilton$-$connected$  if for
any pair of distinct vertices $u$ and $v$, there is a Hamilton $u$-$v$ path, namely, there is
a $u$-$v$ path of length $n-1$. Note that if the graph $G$ is Hamilton-connected then it is hamiltonian. Also note that a hamiltonian graph may not be Hamilton-connected. For instance the cycle $C_6$ is hamiltonian, but it is easy to see that $C_6$ is not  Hamilton-connected. A graph $G$ of order $n > 2$ is $panconnected$  if for every two vertices $u$ and $v$, there is
a $u$-$v$ path of length $l$ for every integer $l$ with
$d(u,v) \leq l \leq n-1$. It is trivial that if $G$ is a panconnected graph, then it is  Hamilton-connected. Alspach [1] showed that the Johnson graph $J(n,m)$ is Hamilton-connected. This result has been generalized in [7] where it has been shown that the Johnson graphs are panconnected. It is easy to check that the graph $L(n)$ is not panconnected. In fact, it is not hard to  show that for
adjacent vertices $v=[1,12]$ and $w=[2,12]$, there is not a $v$-$w$ path   of length 3 (or 4)  in the graph $L(n)$.    In this paper, we wish to show that the graph $L(n)$ has another interesting properties, that is, if $n\geq 6$, then it is  vertex-pancyclic and Hamilton-connected.

 The group of all permutations of a set $V$ is denoted by  $Sym(V)$  or
just  $Sym(n)$ when $ | V | =n $. A $permutation\  group$ $\Gamma$ on
$V$ is a subgroup of  $Sym(V)$. In this case we say that $\Gamma$ acts
on $V$. If $G$ is a graph with vertex-set $V$, then we can view
each automorphism of $G$ as a permutation of $V$, and so $Aut(G)$ is a
permutation group. Let the group $\Gamma$ act  on $V$. We say that $\Gamma$ is
$transitive$ (or $\Gamma$ acts $transitively$  on $V$)   if there is just
one orbit. This means that given any two element $u$ and $v$ of
$V$, there is an element $ \beta $ of  $\Gamma$ such that  $\beta (u)= v
$.

The graph $G$ is called $vertex$-$transitive$  if  $Aut(G)$
acts transitively on $V(G)$. The action of $Aut(G)$ on
$V(G)$ induces an action on $E(G)$, by the rule,
$\beta\{x,y\}$=\\ $\{\beta(x),\beta(y)\}$,
  $\beta\in Aut(G)$, and
$G$ is called $edge$-$transitive$ if this action is
transitive. The graph $G$ is called $symmetric,$  if  for all
vertices $u, v, x, y,$ of $G$ such that $u$ and $v$ are
adjacent, and $x$ and $y$ are adjacent, there is an automorphism
$\alpha$ such that $\alpha(u)=x$,   and $ \alpha(v)=y$. It is clear
that a connected symmetric graph is vertex-transitive and edge-transitive.
 It is easy to show that the graph $L(n)$ is not a symmetric graph. \

Let $\Gamma$ be any abstract finite group with identity $1$, and
suppose $\Omega$ is a subset of   $\Gamma$  with the
properties:
(i) $x\in \Omega \Longrightarrow x^{-1} \in \Omega$,   $ \ (ii)
 \ 1\notin \Omega $. \newline
The $Cayley\  graph$  $G=G (\Gamma; \Omega )$ is the (simple)
graph whose vertex-set and edge-set are defined as follows:
$V(G) = \Gamma $,  $  E(G)=\{\{g,h\}\mid g^{-1}h\in \Omega \}$. As we have already stated, the graph $L(n)$ is a Cayley graph if and only if $n$ is a power of a prime integer [8].

\section{Main results}
\begin{defin} Let $n\geq 3$ be an integer and $[n]=\{1,2,...,n\}$.  The graph $ B(n)$ is a
 graph with the vertex set  $V=\{v \  |  \ v \subset [n] ,  | v |  \in \{ 1,2  \} \} $ and the
edge set $ E= \{ \{ v , w \} \  | \  v , w \in V , v \subset w $ or $ w \subset v \} $. We denote by $L(n)$ the line graph of $B(n)$, that is, $L(n)=L(B(n))$.
\end{defin}
From Definition 2.1, it follows that every vertex of the graph $L(n)$ is of the form $\{\{i\},\{i,j\}  \}$, where $i,j \in [n]$, and $i\neq j$. In this paper, we denote $\{\{i\},\{i,j\}  \}$ by $[i,ij]$. Hence $L(n)$ is the graph with the vertex-set $V(L(n))=V=\{[i,ij], i,j \in [n], i\neq j  \}$, in which two vertices $[i,ij]$ and $[r,rs]$ are adjacent if and only if $i=r$ or $\{i,j \}=\{r,s\}$.
Thus,  if $[i,ij]  $ is a vertex of $ L(n) $, then
$$  N([i,ij])  = \{ [i,ik] \  | \  k \in [n], \  k \neq j,i  \} \cup \{ [j,ij] \},  $$
hence,  $  deg ([i,ij]) = n-1$.  Therefore  $L(n)$ is a regular graph of valency $n-1$.   In fact, $L(n)$ is a vertex-transitive graph [6,8].   By an easy argument, we can show that the graph $L(n)$ is a connected graph with diameter 3. Also, its girth is 3 and hence it is not a bipartite graph.   Figure 1, shows $L(4)$ in the plane.

\definecolor{qqqqff}{rgb}{0.,0.,1.}
\begin{tikzpicture}[line cap=round,line join=round,>=triangle 45,x=1.4cm,y=1.1cm]
\clip(-0.7,-1.) rectangle (7.12,5.5);
\draw (-0.02,0.08)-- (3.08,2.);
\draw (6.56,-0.02)-- (3.08,2.);
\draw (-0.02,0.08)-- (6.56,-0.02);
\draw (-0.12,5.2)-- (2.98,3.3);
\draw (2.98,3.3)-- (6.34,5.18);
\draw (6.34,5.18)-- (-0.12,5.2);
\draw (0.98,3.08)-- (2.12,2.6);
\draw (2.12,2.6)-- (1.,1.8);
\draw (0.98,3.08)-- (1.,1.8);
\draw (4.28,2.62)-- (5.52,3.12);
\draw (4.28,2.62)-- (5.52,2.02);
\draw (5.52,2.02)-- (5.52,3.12);
\draw (-0.79,4.8) node[anchor=north west] {[1,12]};
\draw (5.88,4.92) node[anchor=north west] {[1,13]};
\draw (2.72,4.38) node[anchor=north west] {[1,14]};
\draw (1.04,3.7) node[anchor=north west] {[2,12]};
\draw (-0.21,2.) node[anchor=north west] {[2,23]};
\draw (1.62,2.4) node[anchor=north west] {[2,24]};
\draw (3.76,2.34) node[anchor=north west] {[4,24]};
\draw (5.72,3.26) node[anchor=north west] {[4,41]};
\draw (5.78,2.22) node[anchor=north west] {[4,43]};
\draw (5.96,-0.36) node[anchor=north west] {[3,34]};
\draw (-0.26,-0.22) node[anchor=north west] {[3,23]};
\draw (2.9,1.62) node[anchor=north west] {[3,31]};
\draw (2.12,2.6)-- (4.28,2.62);
\draw (5.52,3.12)-- (2.98,3.3);
\draw (-0.12,5.2)-- (0.98,3.08);
\draw (5.52,2.02)-- (6.56,-0.02);
\draw (1.,1.8)-- (-0.02,0.08);
\draw (6.34,5.18)-- (3.08,2.);
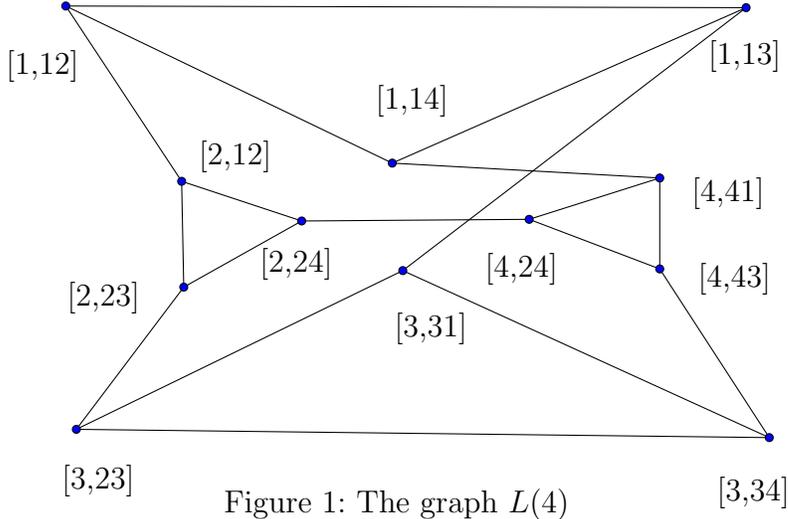
\draw (1.28,-0.5) node[anchor=north west] {Figure 1: The graph $L(4)$};
\begin{scriptsize}
\draw [fill=qqqqff] (2.98,3.3) circle (1.5pt);
\draw [fill=qqqqff] (-0.12,5.2) circle (1.5pt);
\draw [fill=qqqqff] (6.34,5.18) circle (1.5pt);
\draw [fill=qqqqff] (3.08,2.) circle (1.5pt);
\draw [fill=qqqqff] (-0.02,0.08) circle (1.5pt);
\draw [fill=qqqqff] (6.56,-0.02) circle (1.5pt);
\draw [fill=qqqqff] (2.12,2.6) circle (1.5pt);
\draw [fill=qqqqff] (0.98,3.08) circle (1.5pt);
\draw [fill=qqqqff] (1.,1.8) circle (1.5pt);
\draw [fill=qqqqff] (4.28,2.62) circle (1.5pt);
\draw [fill=qqqqff] (5.52,3.12) circle (1.5pt);
\draw [fill=qqqqff] (5.52,2.02) circle (1.5pt);
\end{scriptsize}
\end{tikzpicture} \\

Amongst the various properties of the graph $ L (n)$,
  we interested in  its cycle structure and Hamilton-connectivity. Since for $n=3$ the graph $L(n)$ is isomorphic with $C_6$, and the structure of this graph is simple,  in the sequel we let $n\geq 4$. As we can see in Figure 1,  there are 4 cliques of order 3 in $L(4)$ such that they construct a partition for the vertex-set of this graph. We can easily check that  by these cliques we can construct a $u$-$v$ paths of length $m$ if $5\leq m \leq 11$ for any two  vertices $u$ and $v$ in the graph $L(4)$.  Let $u=[1,12]$ and $v=[2,12]$.  It is easy to check that there is no any $u$-$v$ path of length 4 (or 3)  in $L(4)$. Hence $L(4)$ is not a panconnected graph.
 In the graph $L(n)$ the subgraph induced by the set $C_i=\{[i,ij] | \ j \in [n], j \neq i \}$, $1 \leq i \leq n$,  is a clique of order $n-1$.     It is clear that $P=\{ C_i | \  1 \leq i \leq n \}$ is a partition for the vertex-set of $L(n)$. Note that for any two cliques $C_i$ and $C_j$ there is a unique  pair of adjacent vertices $(u,v)$ such that $u \in C_i$, $v \in C_j$ (indeed $u=[i,ij]$ and $v=[j,ij]$). Moreover if $v$ is a  vertex in the clique $C_i$, then there is exactly one vertex $w$ in $V(L(n))-C_i$ such that $w$ is adjacent to $v$. In fact if $v=[i,ij] \in C_i$, then $w=[j,ij]$ is the unique vertex which is  not in $C_i$ ($w \in C_j$) and adjacent to $v$.

\begin{lem} Let $n \geq 4$ be an integer and $P=\{ C_i | \  1 \leq i \leq n \}$ be the clique partition of
the vertex-set of the graph $L(n)$ where $C_i=\{[i,ij] | \ j \in [n], j \neq i \}$. Let $C \in P$ and let $v,w$ be two vertices in $C$. Let $t$ be an integer
such that,  $2 \leq t \leq n-1$. Then there is a t-subset $P_1=\{A_1,...A_t  \} \subseteq P-\{ C \}$,
and a $v$-$w$ path   $Q:v,v_1,w_1,v_2,w_2,...,v_t,w_t,w$ of length  $2t+1$ such that $v_i,w_i \in A_i$.

\end{lem}

\begin{proof}
 We prove the lemma by induction on $t$. Let $t=2$.  Since $L(n)$ is a vertex-transitive graph, without loss of
generality, we can assume that $v=[1,12]$ (and hence $C=C_1)$, $w=[1,1i]$, $i \neq 1,2$. Consider the vertices
$v_1=[2,12] \in  C_2$ and $w_2=[i,1i] \in C_i$. Thus   there is a unique  adjacent pair $(w_1,v_2)$ such that $w_1
\in C_2$ and $v_2 \in C_i$ ($w_1=[2,2i], v_2=[i,2i] $). Now, it is
clear that the path $Q: v,v_1,w_1,v_2,w_2,w$ is a desired path of length $ 5=2t+1$. Now let $2 \leq m < t$ and the claim is true for $m$.  By induction hypothesis there are cliques $A_1,A_2,...,A_{t-1} \subset P-C_1$ and a $v$-$w$ path   $Q:v,v_1,w_1,v_2,w_2,...,v_{t-1},w_{t-1},w$ of length  $2t-1$ such that $v_i,w_i \in A_i$. Let $A_t \in P- \{C_1,A_1,A_2,...,A_{t-1} \}$ be an arbitrary
clique. There is a unique adjacent pair  $(u_{t-2},z_{t-1})$ such that $u_{t-2} \in A_{t-2}$ and $z_{t-1} \in A_t$. Note that $u_{t-2} \neq v_{t-2},w_{t-2}$. Moreover there is a unique adjacent pair $(u_{t-1},z_{t})$ such that $u_{t-1} \in A_{t}$ and $z_{t} \in A_{t-1}$. Again Note that $z_{t} \neq v_{t-1},w_{t-1}$. Now it is clear that the $v$-$w$ path
 $Q_0:v,v_1,w_1,v_2,w_2,...,v_{t-2},u_{t-2},z_{t-1},u_{t-1},z_t,w_{t-1},w$ of length $2t+1$ is a desired path if we rename $A_{t-1}$ by $A_t$ and $A_t$ by $A_{t-1}$.

\end{proof}  
A graph $G$ of order $n$ is $k$-$pancyclic$ ($k \leq n$) if it contains cycles of every
length from $k$ to $n$ inclusive, and $G$ is $pancyclic$  if it is $g$-pancyclic, where
$g = g(G)$ is the girth of $G$. A graph is of $pancyclicity$ if it is pancyclic.  Not that if $G$  is pancyclic, then $G$ is hamiltonian. The
pancyclicity is an important property to
determine if a topology of a network is suitable for some applications where
mapping cycles of any length into the topology of the network is required [10,11].
The concept of pancyclicity, proposed first by Bondy [5], has been
extended to vertex-pancyclicity  and  edge-pancyclicity  [1]. A graph $G$
of order $n$ is $vertex$-$pancyclic$ (resp. $edge$-$pancyclic$) if any vertex (resp. edge)
lies on cycles of every length from $g(G)$ to $n$ inclusive. Obviously, an edge pancyclic
graph is certainly vertex-pancyclic.

It is easy to check that $L(4)$ has no  cycles of lengths 4 and 5. Also $L(5)$ has no  cycle of length 5.  But by lemma 2.2, we can show that if $n\geq 6$, then $L(n)$ is a vertex-pancyclic graph. Moreover,  we can check that $L(n)$ is not edge-pancyclic, since the edge $e=\{ [1,12], [2,12] \}$ does not lie on a 3-cycle (or 4-cycle). In fact we have the following result.

\begin{thm} Let $n\geq 6$ be an integer. Then $L(n)$ is a vertex-pancyclic graph.

\end{thm}

\begin{proof}
 Let $m$ be an integer such that $n(n-1) \geq m \geq 3$. Let $v$ be an arbitrary vertex in the graph $L(n)$. Let $C$ be the $(n-1)$-clique containing $v$ in $L(n)$. Hence if $3 \leq m \leq 5$, then there is a   cycle $C_m$ containing $v$ in $C$. Now  assume that $6 \leq m \leq 3(n-1)$.  Let $w \in C$ be such that $w \neq v$. By   Lemma 2.2, there are two $(n-1)$-cliques $A_1$ and $A_2$  and a path $Q: v,v_1,w_1,v_2,w_2,w$ of length 5 such that $v_i,w_i \in A_i$ for each $i \in \{  1,2\}$. Let $ 1 \leq t \leq n-3$. It is clear that we can insert $t$ vertices of the clique $C$ between $v$ and $w$. Also, we can insert $t$ vertices of the clique $A_i$, $i \in \{ 1,2 \}$ between $v_i$ and $w_i$. Now it is clear how we can construct an $m$-cycle containing the vertex $v$ in the graph $L(n)$. Now let $m > 3(n-1)$. 
 By Lemma 2.2, there are cliques $A_1,...,A_{(n-1)}$ and the path  $Q:v,v_1,w_1,v_2,w_2,...,v_{(n-1)},w_{(n-1)},w$ in $L(n)$ such that $v_i,w_i \in A_i$, $1 \leq i \leq n-1$. Note that $3(n-1) > 2n$.  By inserting adequate number of vertices of each $A_i$ (and $C$) between each pair of vertices $v_i$ and $w_i$ ($v$ and $w$), we can construct an
 $m$-cycle containing the vertex $v$ in the graph $L(n)$.

\end{proof}

We now want  to show that the graph $L(n)$ is a  Hamilton-connected graph.
\begin{thm}
Let $n \geq 4$ be an integer. Then $L(n)$ is a Hamilton-connected graph.
\end{thm}

\begin{proof} It is easy to  check that the assertion of the theorem is true for the case $n=4$, hence in the sequel we assume that $n \geq 5$.
Let $v,w$ be two vertices in $L(n)$. We show that there is a hamiltotian $v$-$w$ path in $L(n)$.  Let $P=\{ C_i | \  1 \leq i \leq n \}$ be the clique partition of
the vertex-set of the graph $L(n)$ where $C_i=\{[i,ij] | \ j \in [n], j \neq i \}$. There are two cases, namely, \\ (i) $v$ and $w$ are in the same $(n-1)$-clique in the graph $L(n)$, or \\ (ii) $v$ and $w$ are in distinct $(n-1)$-cliques.\newline
 (i) Assume that $A_1 \in P$ and $v,w \in A_1$. Let  $w_1 \in A_1$ be such that $w_1 \neq v,w$.
Therefore by Lemma 2.2,    for subset $P_1=\{A_2,...A_n  \} \subseteq P-\{ A_1 \}$,
there is  a $v$-$w_1$ path   $Q_0:v=v_1,v_2,w_2,...,v_n,w_n,w_1$ of length  $2n-1$ such that $v_i,w_i \in A_i$. Hence  $Q_1:v=v_1,v_2,w_2,...,v_n,w_n,w_1,w$ is a $v$-$w$ path of length $2n$ in $L(n)$. For each $i, 2\leq i \leq n$ we  insert all other vertices in the clique $A_i$ between $v_i$ and $w_i$ and again obtain a $v$-$w$ path of greater length. Now by inserting all the vertices in the set $A_1-\{v,w,w_1  \}$ between the vertices $w_1$ and $w$, we obtain a hamitonian $v$-$w$ path in the graph $L(n)$.\newline
 (ii) We now assume that $v$ and $w$ are in distinct $(n-1)$-cliques. Let $A_1,A_n \in P$ are such that $v \in A_1$ and $w \in A_n$. We know that for each vertex $x$ in $L(n)$  there is a unique
  $(n-1)$-clique $C_x$   in the graph $L(n)$ such that $x\notin C_x $ 
  but  $x$ is adjacent to a unique vertex in $C_x$.     
        Let $P-\{A_1,A_n \}=\{A_2,...,A_{n-1} \}$ be such that $A_2 \neq C_v$ and $A_{n-1} \neq C_w$. Therefore $v$ is adjacent to no vertex in $A_2$ and $w$ is adjacent to no vertex in $A_{n-1}$.  On the other hand,  for each $i$, $1 \leq i \leq n-1$ there is a unique adjacent pair  $w_i, v_{i+1}$ such that $w_i \in A_i$ and $ v_{i+1} \in A_{i+1}$. Note that $w_1 \neq v$ and $v_n \neq w$. Hence the path $Q_0: v=v_1,w_1,v_2,w_2,...v_n,w_n=w$ is a $v$-$w$ path in the graph $L(n)$. For each $i, 1\leq i \leq n$, if  we  insert between $v_i$ and $w_i$ all other vertices in the clique $A_i$,   we obtain a hamitonian $v$-$w$ path in the graph $L(n)$.

\end{proof}

\end{document}